\documentclass{gtart}

\def\ifplaintex{\expandafter\ifx\csname documentclass\endcsname\relax}

\ifplaintex 
\else
\fi

\expandafter\ifx\csname beginpicture\endcsname\relax
\expandafter\ifx\csname documentclass\endcsname\relax
\input pictex \else
\input prepictex \input pictex \input postpictex \fi\fi

\def\gt{{\mathsurround=0pt\it $\cal G\mskip-2mu$eometry \&\ 
$\cal T\!\!$opology}}        

\def\gtp{{\mathsurround=0pt\it $\cal G\mskip-2mu$eometry \&\ 
$\cal T\!\!$opology $\cal P\!$ublications}}  

\def\volumenumber#1{\def\thevolumenumber{#1}}
\def\papernumber#1{\def\thepapernumber{#1}}
\def\volumeyear#1{\def\thevolumeyear{#1}}

\def\pagenumbers#1#2{\def\startpage{#1}\def\finishpage{#2}}
\def\published#1{\def\publishdate{#1}}
\def\proposed#1{\def\theproposer{#1}}
\def\seconded#1{\def\theseconders{#1}}
\def\received#1{\def\receiveddate{#1}}
\def\revised#1{\def\reviseddate{#1}}
\def\accepted#1{\def\accepteddate{#1}}

\let\\\par
\let\thevolumenumber\relax\let\thepapernumber\relax
\let\thevolumeyear\relax\let\thesamplenumber\relax\let\startpage\relax
\let\finishpage\relax\let\publishdate\relax\let\receiveddate\relax
\let\reviseddate\relax\let\accepteddate\relax\let\theasciititle\relax
\let\theasciiauthors\relax
\let\theasciiabstract\relax
\let\theasciiemail\relax\let\theshortauthors\relax\let\theshorttitle\relax

\long\def\maketitlep{   

\count0=\startpage

\gt\hfill      
\beginpicture
\setcoordinatesystem units <0.33truein, 0.33truein> point at 2.2 0.9
\setplotsymbol ({$\cal G$})
\plotsymbolspacing=9truept
\circulararc 315 degrees from 0 1 center at 0 0
\setplotsymbol ({$\cal T$})
\circulararc 315 degrees from 1 -1 center at 1 0
\endpicture
\break
{\small\ifx\thesamplenumber\relax 
Volume \else Sample
\fi\thevolumenumber\ (\thevolumeyear)
\startpage--\finishpage\nl
Published: \publishdate}
\vglue 0.5truein plus 0.4fil minus 0.1truein

{\parskip=0pt\leftskip 0pt plus 1fil\def\\{\par\smallskip}{\ifplaintex\large
\else\Large\fi\bf\thetitle}\par\medskip}   

\vglue 0pt plus 0.1fil 

{\parskip=0pt\leftskip 0pt plus 1fil\def\\{\par}{\sc\theauthors}
\par\medskip}

\vglue 0pt plus 0.1fil 

{\small\parskip=0pt\let\newline\\
{\leftskip 0pt plus 1fil\def\\{\par}{\sl\theaddress}\par}
\expandafter\ifx\theemail\relax    
\relax\else\vglue 5pt plus 0.02fil minus 2pt\def\\{\stdspace{\rm 
and}\stdspace} 
\cl{Email:\stdspace\tt\theemail}\fi
\ifx\theurl\relax                  
\relax\else\vglue 5pt plus 0.02fil minus 2pt\def\\{\stdspace{\rm 
and}\stdspace}
\cl{URL:\stdspace\tt\theurl}\fi\par}

\vglue 7pt plus 0.3fil minus 3pt

{\bf Abstract}
\vglue 5pt plus 0.1fil minus 2pt

\theabstract

\vglue 7pt plus 0.3fil minus 3pt

{\bf AMS Classification numbers}\quad Primary:\quad \theprimaryclass

Secondary:\quad \thesecondaryclass

\vglue 5pt plus 0.3fil minus 2pt

{\bf Keywords:}\quad \thekeywords

\vglue 10pt plus 0.5fil minus 5pt

{\small  Proposed: \theproposer\hfill Received: \receiveddate\nl
Seconded: \theseconders\hfill 
\ifx\reviseddate\relax                         
Accepted: \accepteddate                        
\else
Revised: \reviseddate                          
\fi}
\eject
}       

\let\maketitlepage\maketitlep
\let\maketitle\maketitlepage

\font\phead=cmsl9 scaled 950
\font=cmsl9 scaled 1050
\font\pnum=cmbx10 scaled 913
\font\lnum=cmbx10 
\font\pfoot=cmsl9 scaled 950
\font=cmsl9 scaled 1050
\ifplaintex
\headline{\vbox to 0pt{\vskip -4.5mm\line{\small\phead\ifnum
\count0=\startpage ISSN 1364-0380 (on line)
1465-3060 (printed) \hfill {\pnum\folio}\else\ifodd\count0\def\\{ }%
\ifx\theshorttitle\relax\thetitle\else\theshorttitle\fi\hfill{\pnum\folio}
\else\def\\{ and }{\pnum\folio}\hfill\ifx\theshortauthors\relax\theauthors
\else\theshortauthors\fi\fi\fi}\vss}}
\footline{\vbox to 0pt{\vglue 0mm\line{\small\pfoot\ifnum\count0=\startpage
\copyright\ \gtp\hfill\else
\gt, Volume \thevolumenumber\ (\thevolumeyear)\hfill\fi}\vss
}}
\else
\makeatletter
\def\@oddhead{{\small\lhead\ifnum\count0=\startpage ISSN 1364-0380 (on line)
1465-3060 (printed) \hfill {\lnum\number\count0}\else\ifodd\count0
\def\\{ }\ifx\theshorttitle\relax \thetitle \else\theshorttitle\fi\hfill
{\lnum\number\count0}\else\def\\{ and }{\lnum\number\count0}
\hfill\ifx\theshortauthors\relax 
\theauthors\else\theshortauthors\fi\fi\fi}}\def\@evenhead{\@oddhead}
\def\@oddfoot{\small\lfoot\ifnum\count0=\startpage\copyright\ \gtp\hfill\else
\gt, Volume \thevolumenumber\ (\thevolumeyear)\hfill\fi}
\def\@evenfoot{\@oddfoot}
\makeatother
\fi

\newwrite\gtoutfile
\long\gdef\makeheadfile{  
{\def\\{, }\def\s{ }
\immediate\openout\gtoutfile head.xxx
\immediate\write\gtoutfile{To: math@arxiv.org}
\immediate\write\gtoutfile{Subject: put or rep NNNNN:pppp}
\immediate\write\gtoutfile{--text follows this line--}
\immediate\write\gtoutfile{Proxy-for: \ifx\theasciiauthors\relax
\theauthors\else\theasciiauthors\fi\s<\ifx\theasciiemail\relax\theemail\else\theasciiemail\fi>}
\immediate\write\gtoutfile{\noexpand\\}
\immediate\write\gtoutfile{Authors: \ifx\theasciiauthors\relax
\theauthors\else\theasciiauthors\fi}
\immediate\write\gtoutfile{Title: \ifx\theasciititle\relax
\thetitle\else\theasciititle\fi}
\immediate\write\gtoutfile{Subj-class: GT or SG or MG etc}
\immediate\write\gtoutfile{MSC-class: \theprimaryclass\ifx\thesecondaryclass\relax\else, \thesecondaryclass\fi}
\immediate\write\gtoutfile{Journal-ref: Geom. Topol. \thevolumenumber
(\thevolumeyear) \startpage-\finishpage}
\immediate\write\gtoutfile{Comments: Published in Geometry and Topology at}
\immediate\write\gtoutfile{\s\s http://www.maths.warwick.ac.uk/gt/GTVol\thevolumenumber/paper\thepapernumber.abs.html}
\immediate\write\gtoutfile{\noexpand\\}
\immediate\write\gtoutfile{}
\ifx\theasciiabstract\relax
\immediate\write\gtoutfile{\theabstract}\else
\immediate\write\gtoutfile{\theasciiabstract}\fi
\immediate\write\gtoutfile{}
\immediate\write\gtoutfile{\noexpand\\}
\immediate\write\gtoutfile{}
\immediate\closeout\gtoutfile}}  

\def\maketitlepage{\maketitlep\makeheadfile}
\let\maketitle\maketitlepage

\def\ifplaintex{\expandafter\ifx\csname documentclass\endcsname\relax}

\ifplaintex 
\else
\fi

\expandafter\ifx\csname beginpicture\endcsname\relax
\expandafter\ifx\csname documentclass\endcsname\relax
\input pictex \else
\input prepictex \input pictex \input postpictex \fi\fi

\def\gt{{\mathsurround=0pt\it $\cal G\mskip-2mu$eometry \&\ 
$\cal T\!\!$opology}}        

\def\gtp{{\mathsurround=0pt\it $\cal G\mskip-2mu$eometry \&\ 
$\cal T\!\!$opology $\cal P\!$ublications}}  

\def\volumenumber#1{\def\thevolumenumber{#1}}
\def\papernumber#1{\def\thepapernumber{#1}}
\def\volumeyear#1{\def\thevolumeyear{#1}}

\def\pagenumbers#1#2{\def\startpage{#1}\def\finishpage{#2}}
\def\published#1{\def\publishdate{#1}}
\def\proposed#1{\def\theproposer{#1}}
\def\seconded#1{\def\theseconders{#1}}
\def\received#1{\def\receiveddate{#1}}
\def\revised#1{\def\reviseddate{#1}}
\def\accepted#1{\def\accepteddate{#1}}

\let\\\par
\let\thevolumenumber\relax\let\thepapernumber\relax
\let\thevolumeyear\relax\let\thesamplenumber\relax\let\startpage\relax
\let\finishpage\relax\let\publishdate\relax\let\receiveddate\relax
\let\reviseddate\relax\let\accepteddate\relax\let\theasciititle\relax
\let\theasciiauthors\relax
\let\theasciiabstract\relax
\let\theasciiemail\relax\let\theshortauthors\relax\let\theshorttitle\relax

\long\def\maketitlep{   

\count0=\startpage

\gt\hfill      
\beginpicture
\setcoordinatesystem units <0.33truein, 0.33truein> point at 2.2 0.9
\setplotsymbol ({$\cal G$})
\plotsymbolspacing=9truept
\circulararc 315 degrees from 0 1 center at 0 0
\setplotsymbol ({$\cal T$})
\circulararc 315 degrees from 1 -1 center at 1 0
\endpicture
\break
{\small\ifx\thesamplenumber\relax 
Volume \else Sample
\fi\thevolumenumber\ (\thevolumeyear)
\startpage--\finishpage\nl
Published: \publishdate}
\vglue 0.5truein plus 0.4fil minus 0.1truein

{\parskip=0pt\leftskip 0pt plus 1fil\def\\{\par\smallskip}{\ifplaintex\large
\else\Large\fi\bf\thetitle}\par\medskip}   

\vglue 0pt plus 0.1fil 

{\parskip=0pt\leftskip 0pt plus 1fil\def\\{\par}{\sc\theauthors}
\par\medskip}

\vglue 0pt plus 0.1fil 

{\small\parskip=0pt\let\newline\\
{\leftskip 0pt plus 1fil\def\\{\par}{\sl\theaddress}\par}
\expandafter\ifx\theemail\relax    
\relax\else\vglue 5pt plus 0.02fil minus 2pt\def\\{\stdspace{\rm 
and}\stdspace} 
\cl{Email:\stdspace\tt\theemail}\fi
\ifx\theurl\relax                  
\relax\else\vglue 5pt plus 0.02fil minus 2pt\def\\{\stdspace{\rm 
and}\stdspace}
\cl{URL:\stdspace\tt\theurl}\fi\par}

\vglue 7pt plus 0.3fil minus 3pt

{\bf Abstract}
\vglue 5pt plus 0.1fil minus 2pt

\theabstract

\vglue 7pt plus 0.3fil minus 3pt

{\bf AMS Classification numbers}\quad Primary:\quad \theprimaryclass

Secondary:\quad \thesecondaryclass

\vglue 5pt plus 0.3fil minus 2pt

{\bf Keywords:}\quad \thekeywords

\vglue 10pt plus 0.5fil minus 5pt

{\small  Proposed: \theproposer\hfill Received: \receiveddate\nl
Seconded: \theseconders\hfill 
\ifx\reviseddate\relax                         
Accepted: \accepteddate                        
\else
Revised: \reviseddate                          
\fi}
\eject
}       

\let\maketitlepage\maketitlep
\let\maketitle\maketitlepage

\font\phead=cmsl9 scaled 950
\font=cmsl9 scaled 1050
\font\pnum=cmbx10 scaled 913
\font\lnum=cmbx10 
\font\pfoot=cmsl9 scaled 950
\font=cmsl9 scaled 1050
\ifplaintex
\headline{\vbox to 0pt{\vskip -4.5mm\line{\small\phead\ifnum
\count0=\startpage ISSN 1364-0380 (on line)
1465-3060 (printed) \hfill {\pnum\folio}\else\ifodd\count0\def\\{ }%
\ifx\theshorttitle\relax\thetitle\else\theshorttitle\fi\hfill{\pnum\folio}
\else\def\\{ and }{\pnum\folio}\hfill\ifx\theshortauthors\relax\theauthors
\else\theshortauthors\fi\fi\fi}\vss}}
\footline{\vbox to 0pt{\vglue 0mm\line{\small\pfoot\ifnum\count0=\startpage
\copyright\ \gtp\hfill\else
\gt, Volume \thevolumenumber\ (\thevolumeyear)\hfill\fi}\vss
}}
\else
\makeatletter
\def\@oddhead{{\small\lhead\ifnum\count0=\startpage ISSN 1364-0380 (on line)
1465-3060 (printed) \hfill {\lnum\number\count0}\else\ifodd\count0
\def\\{ }\ifx\theshorttitle\relax \thetitle \else\theshorttitle\fi\hfill
{\lnum\number\count0}\else\def\\{ and }{\lnum\number\count0}
\hfill\ifx\theshortauthors\relax 
\theauthors\else\theshortauthors\fi\fi\fi}}\def\@evenhead{\@oddhead}
\def\@oddfoot{\small\lfoot\ifnum\count0=\startpage\copyright\ \gtp\hfill\else
\gt, Volume \thevolumenumber\ (\thevolumeyear)\hfill\fi}
\def\@evenfoot{\@oddfoot}
\makeatother
\fi

\newwrite\gtoutfile
\long\gdef\makeheadfile{  
{\def\\{, }\def\s{ }
\immediate\openout\gtoutfile head.xxx
\immediate\write\gtoutfile{To: math@arxiv.org}
\immediate\write\gtoutfile{Subject: put or rep NNNNN:pppp}
\immediate\write\gtoutfile{--text follows this line--}
\immediate\write\gtoutfile{Proxy-for: \ifx\theasciiauthors\relax
\theauthors\else\theasciiauthors\fi\s<\ifx\theasciiemail\relax\theemail\else\theasciiemail\fi>}
\immediate\write\gtoutfile{\noexpand\\}
\immediate\write\gtoutfile{Authors: \ifx\theasciiauthors\relax
\theauthors\else\theasciiauthors\fi}
\immediate\write\gtoutfile{Title: \ifx\theasciititle\relax
\thetitle\else\theasciititle\fi}
\immediate\write\gtoutfile{Subj-class: GT or SG or MG etc}
\immediate\write\gtoutfile{MSC-class: \theprimaryclass\ifx\thesecondaryclass\relax\else, \thesecondaryclass\fi}
\immediate\write\gtoutfile{Journal-ref: Geom. Topol. \thevolumenumber
(\thevolumeyear) \startpage-\finishpage}
\immediate\write\gtoutfile{Comments: Published in Geometry and Topology at}
\immediate\write\gtoutfile{\s\s http://www.maths.warwick.ac.uk/gt/GTVol\thevolumenumber/paper\thepapernumber.abs.html}
\immediate\write\gtoutfile{\noexpand\\}
\immediate\write\gtoutfile{}
\ifx\theasciiabstract\relax
\immediate\write\gtoutfile{\theabstract}\else
\immediate\write\gtoutfile{\theasciiabstract}\fi
\immediate\write\gtoutfile{}
\immediate\write\gtoutfile{\noexpand\\}
\immediate\write\gtoutfile{}
\immediate\closeout\gtoutfile}}  

\def\maketitlepage{\maketitlep\makeheadfile}
\let\maketitle\maketitlepage

\volumenumber{4}\papernumber{15}\volumeyear{2000}
\pagenumbers{431}{449}
\proposed{Cameron Gordon}
\seconded{David Gabai, Robion Kirby}
\received{20 February 1999}
\revised{29 May 2000}
\accepted{11 November 2000}
\published{14 November 2000}

\usepackage{amssymb, amsmath, graphicx, psfrag, labelfig} 

\newtheorem{theorem}{Theorem}[section]    
\newtheorem{lemma}[theorem]{Lemma}

\newcommand{\ZZ}{\mathbb{Z}}

\newcommand{\qp}{\mathbb{QP}^1}
\newcommand{\fp}{\mathbb{F}_p\mathbb{P}^1}

\def\del{$\partial$}
\def\a{$\alpha$}
\def\b{$\beta$}
\def\g{$\gamma$}

\def\L{$\mathcal{L}$\ }

\def\min{\mbox{\rm{min}}}

\def\mod{\mbox{mod}}

\def\H{{\mathbb H}^3}

\begin{document}

\title{Bounds on exceptional Dehn filling}                    
\authors{Ian Agol}                  
\address{Department of Mathematics\\University of Melbourne\\Parkville, VIC 3052\\Australia}

\email{agol@ms.unimelb.edu.au}                     
\url{http://math.ucdavis.edu/\char126 iagol}                       

\begin{abstract}   
We show that for a hyperbolic knot complement, all but at most 12 Dehn
fillings are irreducible with infinite word-hyperbolic fundamental
group.
\end{abstract}

\primaryclass{57M50, 57M27}                
\secondaryclass{57M25, 57S25}              
\keywords{Hyperbolic, Dehn filling, word-hyperbolic}                    

\maketitlepage

\section{Introduction}
Thurston demonstrated that if one has a hyperbolic knot complement,
all but finitely many Dehn fillings give hyperbolic manifolds \cite{Th}. The 
example with the largest known number of non-hyperbolic Dehn 
fillings is the figure-eight knot complement, which has 10 
fillings which are not hyperbolic. It is conjectured that this
is the maximal number that can occur. Call a manifold {\it hyperbolike}
if it is irreducible with infinite word-hyperbolic fundamental 
group (this is stronger than Gordon's definition \cite{G}). For
example, manifolds with a Riemannian metric of negative sectional
curvature are hyperbolike.
Call a Dehn filling {\it exceptional} if it is not hyperbolike.
We will consider the more amenable problem of determining
the number of exceptional Dehn fillings on a knot complement. 
The geometrization conjecture would imply that hyperbolike manifolds
are hyperbolic.
Bleiler and Hodgson \cite{BH} showed that there are at most 24 exceptional 
Dehn fillings, using Gromov and Thurston's $2\pi$--theorem and
estimates on cusp size due to Colin Adams\cite {A}. We will make an
improvement on the $2\pi$--theorem, and use improved lower bounds
on cusp size due to Cao and Meyerhoff \cite{CM}, to get an upper bound
of 12 exceptional Dehn fillings. The inspiration for this work
came from discussions with Zheng-Xu He. He has obtained bounds relating
asymptotic crossing number to cusp geometry \cite{He}. He remarked
to me that his estimates could probably be improved, and this
paper gives my attempt at such an improvement. Mark Lackenby 
\cite{L} has
obtained the same improvement of the $2\pi$--theorem. My 
thesis gave implications about Dehn fillings being atoroidal, not 
word-hyperbolic \cite{Ag}. Marc has an improved
version of Gabai's ubiquity theorem which filled in a gap in an
early draft of my thesis, and this is the argument which appears in
this paper.

\section{Definitions}
The notation introduced here will be used throughout the paper. We will
use ${\rm int} X$ to mean the interior of the space $X$, and ${\mathcal N}(X)$
will denote an open regular neighborhood of a subset $X\subset M$. 
$M$ is a hyperbolic $3$--manifold  with a distinguished torus cusp. $M$ has
a compactification to a compact $3$ manifold $\overline{M}$ with torus boundary,
by adding the ends of geodesic rays which remain in the cusp. 
Let $S$ be a surface of finite type ($S$ may have both boundary and punctures), let
 $f\co S\to M$ be
 a
 map such that every puncture maps properly into a cusp. This map might not
necessarily be an embedding or an immersion. Using the terminology
 of 
\cite{O}, $f\co S\to M$ is {\it incompressible} if every simple loop $c$ in $S$ 
 for which $f(c)$ is homotopically trivial in $M$   
bounds a disk in $S$. The simple loop conjecture would imply that such a map $f$ 
is $\pi_1$--injective, but this is not known for general $3$--manifolds \cite{Ga1}. 
Let $U=[0,\infty)\times \mathbf{R}$.
 A boundary compression of $f$ is a proper map $b\co U\to M$ such that there is a map
$b'\co \partial U\to S$ with $f\circ b'=b_{|\partial U}$, $b'(\partial U)$ 
is a proper simple line in $S$ which does not bound a properly embedded
half-plane in $S$. 
 $f$ is {\it \del--incompressible} if it has no boundary compression. $f\co S\to M$ is
{\it essential} if it is incompressible and \del--incompressible. $f\co S\to M$ is 
{\it pleated} 
if the boundary components of $S$ map to geodesics in $M$, and ${\rm int} S$ (the interior
of $S$) is piecewise 
made of triangles which map under $f$ to ideal hyperbolic geodesic triangles in $M$, so 
that the 1--skeleton $\cup\ \partial S$ forms a lamination in
 $S$. A pleated surface has an induced hyperbolic metric, where the lamination 
is geodesic. 

For a cusped hyperbolic $3$--manifold $M$, we may take an embedded neighborhood
$C$ of the cusp  which is a  quotient of an open horoball by the torus group, which
we will call a {\it horocusp}. The closure of $C$ might not be embedded, so by 
 $\partial C$ we will mean the torus obtained as the path closure of the Riemannian manifold
$C$ (not regarded as a subset of $M$). $\partial C$ inherits a euclidean metric from $M$. 
 If $p$ a loop in 
$C$, let $l_C(p)$ denote the length of a euclidean geodesic loop  homotopic to $p$
in $\partial C$. A {\it slope} in $\partial C$ is an equivalence class of embedded loops in
$\partial C$. If $\alpha$ is a slope in $\partial C$, then $M(\alpha)$ denotes the Dehn
filling along that slope, which is a manifold obtained by gluing a solid torus to
$M\backslash C$ so that the loop represented by $\alpha$ bounds a disk in the solid torus.
This is uniquely determined by the slope in $\partial C$. 

A theorem of Gromov \cite{Gr} implies that for a closed manifold $M$,
$\pi_1(M)$ is {\it word-hyperbolic} if for a metric on $M$, $M$ has a
linear isoperimetric inequality.  That is for a metric on $M$, there
is a constant $V$ so that for any map of a disk $d\co D^2\to M$,
$\text{ area}(D)\leq V \text{ length}(\partial D)$ in the induced
metric on $D$. Gromov has shown that for such a manifold, $\pi_1M$ has
no $\ZZ+\ZZ$ subgroup and has a solvable word problem. A theorem of
Bestvina and Mess implies that the universal cover $\tilde M$ has a
compactification to a ball \cite{BM} (if the fundamental group is
infinite and the manifold is irreducible). Thurston's geometrization
conjecture would imply that $M$ has a hyperbolic structure, that is a
Riemannian metric of constant sectional curvature -1.

\section{Essential Surfaces}
In this section, we show how to obtain singular essential surfaces in a knot
complement coming from the ambient manifold. The results are similar to those
of Ulrich Oertel \cite{O}, but we do not worry about embeddedness of the
boundary components.  Marc Lackenby \cite{L} and Zheng-Xu He \cite{He} have also
obtained similar results to the following lemmas.
The idea is that if there is an essential sphere or if the core of
a Dehn filling has finite order in the fundamental group, so that 
some multiple of the core bounds a disk,  then the surface  can be
homotoped so that its intersection with the knot complement
is essential. 

\begin{lemma}[Essential punctured spheres]\label{essential surfaces}
Let $M^3$ be a compact 3--manifold and take a knot $k\subset M$ with 
$N=M\setminus \mathcal N(k)$, such that $\partial N$ is 
incompressible in $N$ and
N is irreducible. Let 
$f\co S\to M$ be a singular map of a sphere or disk. If $S$ is a sphere, then $f$ is a
homotopically non-trivial map into $M$. 
 If $S$ is a disk, then its 
boundary maps to a homotopically non-trivial curve in $\mathcal N (k)$. Then we can find a 
surface $T$ and a mapping  
$g\co T\to M$, with the same properties as above, such that $g$ is transverse to
$\partial N$ and $g^{-1}(N)$ is 
essential in N.
\end{lemma}
\proof
First notice that we may take $f$ transverse to $k$, so that $f^{-1}{\mathcal N (k)}$ is a 
collection of disks in $S$  and an annular neighborhood of
the boundary in the disk case (we will call these {\it dots}). Then we will induct on 
$|f^{-1}{\mathcal N (k)}| =$ the number of dots. 
Suppose there is an essential simple closed curve $c$ in $f^{-1}(N) \equiv \hat{S}$, which bounds 
a disk
$D$ in $N$ (ie, there is a homeomorphism $d'\co c\to \partial D$ and a map $d\co D\to N$ with 
$d_{|\partial D}\circ d'=f$). Then since $S$ is either a disk or
a sphere, $c$ bounds a disk $E$ in $S$ which must
meet $\mathcal N(k)$, since $c$ is essential in $\hat S$. Surger $f\co S\to M$ along $d\co D\to M$.
That is, create a new surface $S'$ by splitting $S$ along $c$, and glue in two copies
of $D$ to the two new boundary components (corresponding to two copies of c) by gluing 
using the homeomorphism $d'\co c\to\partial D$, 
then form a mapping $f'\co S'\to M$ by using $f$ or $d$ on the relevant pieces of $S'$.
One  component of $S'$ may have image in $M$ under $f'$ a homotopically trivial sphere, 
so we get rid of it. In case $S$ is a sphere, there are two 
choices for the disk $E$ bounding $c$ in $S$. At least one choice will surger to
an essential sphere in $M$, so we keep this one. We then have a surface which has fewer dots.
Replace $S$ with this surface, which we will still call $S$, and $f$ with the restriction
of $f'$ to this subsurface, which we will still call $f$.
	
Suppose there is an arc \a\ which is embedded and essential
 in $\hat S$ which bounds
a boundary compression for $\hat S$. That is, there is a map $d'\co \partial U\to \alpha\subset \hat S$
and a map $d\co U\to N$ such that $f\circ d'=d_{|\partial U}$. There are two types of boundary 
compressions:
\begin{enumerate}
\item \a\ connects different dots in $S$, so we push $S$ along 
the boundary compression. That is, we split $S$ along the arc $\overline\alpha$ and make
a new surface $S'$ by gluing two copies of the disk $\overline U$ to the new boundary components
using the homeomorphism $d'\co \partial U\to \alpha$, and identifying the other 
ends of the two copies of $\partial \overline U$ by the identity. Then replace $f$ with $f'\co S'\to M$
by using $f$ or $d$ on the relevant pieces of $S'$. Then, we may expand $\mathcal N (k)$
slightly, and make $f'$ transverse to $\partial N$. 
 This has the effect of turning two dots into one, and decreases
the number of dots. By induction, we may assume there are no such compressions. 

\item \a\ connects the same dot in $S$. Take a maximal collection of disjoint,
non-parallel \del--compres\-sions, and as in case 1, we push $S$ along these \del--compres\-sions
(see the previous case for a more precise description of this push operation), getting
a new map which we will still call $f\co S\to M$. 
Then $f^{-1}(\mathcal N(k))$ is a collection of planar surfaces such that each one 
separates $S$. Take
an innermost planar surface. If there are no dots in the disks it separates off, or
if there are disks whose boundary maps to a homotopically trivial loop in $\mathcal
N(k)$, then 
we can homotope $f$ in a neighborhood of these disks in $S$ into $\mathcal{N}(k)$,
keeping $f$ fixed on  the rest of $S$, since \del $N$ is incompressible in $N$
and $N$ is irreducible,
decreasing the number of dots in $S$. 
Otherwise, one of these disks has boundary which is essential in $\mathcal N(k)$.
We then take T to be this disk adjoined a collar of the outermost curve in $\mathcal{N}(k)$,
and $g=f_{|T}$.\endproof
\end{enumerate}\medskip

The next lemma deals with the case in which one has a singular map of a disk into $M$, 
with boundary mapped into the complement of the knot. Then one can homotope the map
of the disk
to be essential in the knot complement, as long as there are no essential punctured
disks in the knot complement. This will be used later for bounding the area
of such a disk. 

\begin{lemma}[Essential punctured disks] \label{essential disks}
Let $M^3$ be a compact 3--manifold, and take a knot $k\subset M$ with $N=M\setminus \mathcal 
N(k)$, such that $\partial N\subset N$ is incompressible and $N$ is
irreducible. 
Also, assume that there are no maps of disks $a\co A\to M$ transverse to $k$, with 
$a(\partial A) \subset \mathcal N(k)$, and $a^{-1}(N)$ essential
in N. Then, if $f\co D\to M$ is a disk whose boundary is in N, we may homotope $f$ so that
$f^{-1}(N)$ is essential in N. 
\end{lemma}
\begin{proof}
The proof is similar to that of the previous lemma, but we need to observe that
$\pi_2 M = 0$ by the previous lemma, so the disk surgeries can be done by 
homotopies. We may assume that $f_{|f^{-1}(N)}$ is incompressible
and has no \del--compression such that the arc $\alpha$ connnects different dots of 
$f^{-1}(\mathcal N(k))$. 
As before, take a maximal collection
of non-parallel arcs which separate $D$ and bound \del--compressions, and push $f\co D\to M$ along these
\del--compressions to obtain planar surfaces separating $D$
(as in the previous lemma). None of the disks separated
by an innermost planar surface can be essential, since this would contradict our assumption.
Thus, as in the previous lemma, we may homotope 
$f$ on the innermost disks into $\mathcal N (k)$,
decreasing the number of dots. 
\end{proof}

\section{Pleated Surfaces}
The argument in this section is similar to that of Thurston
\cite{Th2}, but 
the hypotheses are slightly different. This result will be used
next section to compare the geometry of the hyperbolic metric on the 
pleated surface to the geometry of the manifold.

\begin{lemma}[Pleated Surfaces] \label{pleated surfaces}
Let $N$ be a  hyperbolic $3$--manifold with a distinguished cusp, let 
$S$ be a surface of finite type with $\chi(S)<0$, and let $f\co S\to N$ be
a singular essential map,
 with cusps of $S$ mapping properly to cusps of $N$, and 
\del $S$ mapping to geodesics in $N$. Then we can find a hyperbolic metric on $S$ and
a map $g\co S\to N$, such that $g$ is pleated, $g_{|{\rm int} S}$ is homotopic to $f_{|{\rm int} S}$, and
$g_{|\partial S}$ is an isometry.\end{lemma}

\begin{proof}Choose an ideal triangulation 
$T$ of ${\rm int} S$, such that no edges connect \del $S$ to itself and
every edge is essential in $S$. Then 
spin the triangles of 
$T$ around \del $S$. We obtain a lamination \L consisting
of $T^{(1)}\cup \partial S$, which is a geodesic lamination, in the sense that it is 
isotopic to a geodesic lamination in any complete hyperbolic structure on $S$ with geodesic
boundary. We may assume $f\co S\to N$ is $C^2$ near \L, since it can be assumed that
the singularities of $f$ are in the interior of $S$, so we can make \L miss the 
singularities.   So each end of a leaf of \L which limits to
\del $S$ must eventually  have curvature close to 0, 
and is therefore a quasi-geodesic in $N$. Its other end
must map properly into a cusp of $N$. Lifting to $\H\cong \tilde{N}$, we see that the endpoints must be 
distinct on \del$\H$, by discreteness. If both ends of a leaf $L$ of \L map into the same cusp
when lifted to $\H$, then $L$ bounds a \del--compressing disk in $\H$, whose end maps into
the same cusp. Pushing down to $N$, we find a \del--compression for $S$, which contradicts 
that the edges of T are essential in $S$ and $S$ is \del--incompressible. In either case, the
endpoints of each leaf lifted to $\H$ map to different points in \del$\H$, so each leaf is
homotopic to a unique geodesic. Therefore, we may homotope $f$ so that $T^{(1)}$ is geodesic in $N$. 
We can homotope $f$ on each  triangle of $T$ to be totally geodesic by homotopy extension
in ${\mathbb H}^3$ and pushing down to $N$, giving a homotopic pleated map  
$g\co {\rm int} S\to N$. A pleated surface has an induced hyperbolic metric, which we give to 
${\rm int} S$. Then we can complete the metric on ${\rm int} S$ to a metric on a 
surface $S'\cong S$. Choose a 
geodesic $\gamma$ in $f(\partial S)$. Then since the ends of leaves of \L are quasi-geodesic,
each end of a geodesic leaf of $g(\mathcal L)$ is in a bounded neighborhood of the end in 
$f(\mathcal L)$. Therefore, the ends of $g(\mathcal L)$ limit to $\gamma$. When part of 
a geodesic of \L wraps closely once about \del $S'$, then its image wraps closely about $\gamma$
in $N$. In the limit, we see that the length of \del $S'$ must be the same as the length of 
$\gamma$. So we may extend $g$ to an isometry $g\co S\to N$.
\end{proof}

\section{Cusp area shrinks}
The next theorem is based on the fact that there are disjointly embedded cusps in 
$S$ which have longer boundary
 than the cusp lengths in the image. This would be easy to show if $S$ were totally geodesic, and
we would get equality. But since $S$ is actually pleated, the folding makes $S$ have longer 
cusp lengths. 

\begin{theorem}[Bounding cusp length by Euler characteristic] \label{bound} Let
$N$ be a\break hyperbolic 3--manifold with a distinguished horocusp $C$. Let 
$S$ be a surface of finite type with no boundary components, and $f\co S\to M$ be  an  essential mapping,
where the cusps of $S$ map into $C$. 
For each puncture $p_i$ of $S$, consider the length of 
the corresponding slope
 $l_C(p_i)$ in $\partial C$. Then $\Sigma_i\ l_C(p_i) \leq 6\  |\chi (S)|$.
\end{theorem}

\begin{proof}
First, homotope $f$ to a pleated map by lemma \ref{pleated surfaces}, which we
will also call $f$. 
$f({\rm int} S)$ is a union of ideal geodesic triangles $T_1$,...,$T_{2|\chi(S)|}$. 
If a corner of $T_j$ is 
 in $C$, the opposite edge of $T_j$ might intersect $C$ in its interior. Lifting $T_j$ to $\tilde{T}_j$
 in ${\tilde{M}} = \mathbb{H}^3$, it looks like Figure \ref{pleated}, with a parabolic limit point
 of $C$ lifted to $\infty$, and $\tilde{C}$ a lift of $C$, in the upper half-space model of
 $\mathbb H^3$. Shrink the cusp $C$ to a cusp $C'$ such that each edge of  $T_j$ intersects
$C$ in no compact intervals. Then $\tilde{C}'$ 
looks like Figure \ref{pleated}.
\begin{figure}[ht!]
\cl{\small
\SetLabels 
\E(.05*.42){$\tilde{C}^{\prime}$}\\
\E(0.47*.68){$T_j$}\\
\E(.13*.06){$\mathbb{H}^3$}\\
\E(.09*.18){$\tilde{C}$}\\
\endSetLabels 
\AffixLabels{\includegraphics[width=.9\hsize]{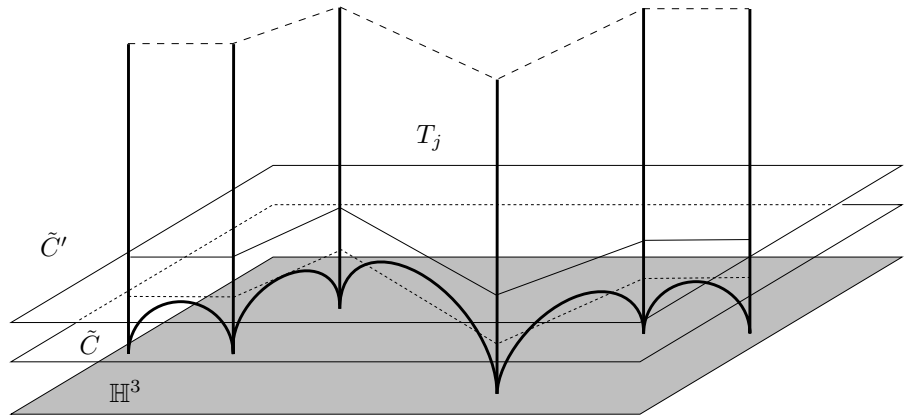}}}
\caption{\label{pleated} Shrinking the horocusp}
\end{figure}
$f^{-1}(C')=H'=\cup H_i'$ consists of disjoint horocusps $H_i'$ in $S$, one for each puncture 
of $S$ which maps into $C$. Let $l_{H'}(p_i) = $ the length of $p_i$ along $\partial H_i'$ in $S$. Then
 $l_{H'}(p_i) \geq l_{C'}(f(p_i))$, with equality iff 
there is no bending along the pleats at $p_i$. Let $d=d(C,C')$. Then choose horocusps $H_i
 \supseteq H_i'$ in S, such that $d(H_i,H_i') = d$. Then $f(H_i) \subseteq C$, since $f$ is 
piecewise an isometry, so it shrinks distances. Suppose $H_i\cap H_j \neq \emptyset$,
 for some $i\neq j$. Then there is a geodesic arc $a$ in $S$ connecting $p_i$ to $p_j$: just
 look at intersecting lifts of $H_i,H_j$ in $\tilde{S} = \mathbb H^2$, and take the 
geodesic $\tilde a$ connecting the centers of 
 $\tilde{H}_i,\tilde{H}_j$ in $\partial{\mathbb H}^2$. See Figure  \ref{horo}.
\begin{figure}[ht!]\small
	\begin{center}
	\psfrag{H}{$\mathbb{H}^2$}
	\psfrag{H1}{$\widetilde{H}_i$}
	\psfrag{H2}{$\widetilde{H}_j$}
	\psfrag{a}{$\tilde{a}$}
	\includegraphics[width=3.2in]{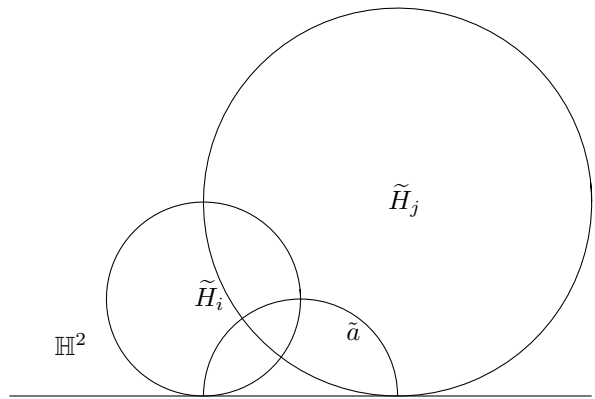}
	\caption{\label{horo} Intersecting horocusps give \del--compression}
	\end{center}
\end{figure}
 $f(H_i), f(H_j) \subseteq C$, so $f(a) \subseteq C$. Then there is a \del--compressing
 disk $D$ for $f(a)$ in $C$. Just cone off $f(a)$ to the end of $C$ by geodesics. That
 is, in the universal cover of $N$, take a cover $\tilde{C}$ of $C$ tangent to $\infty$, and 
$\widetilde{f(a)}$ of $f(a)$. Then over each point of $f(a)$, take a geodesic connecting the point
 to $\infty$. This describes a map of a half plane compressing $\widetilde{f(a)}$. Map
 down to $N$, to get a \del--compression of the arc $a$ in $f$. 

So we have shown that since $S$ is \del--incompressible, $H_i \cap H_j = \emptyset, i\neq j$. 
Thus, we have disjoint horocusps in $S$. 
 $l_H(p_i) = e^d l_{H'}(p_i)$ and $l_C(f(p_i)) = e^d l_{C'}(f(p_i))$. So 
$$l_H(p_i) = e^d l_{H'}(p_i) \geq e^d l_{C'}(f(p_i)) = l_C(f(p_i)).$$ 
A theorem of Boroczky 
\cite{B} implies that ${\rm Area}(H)\leq \frac{3}{\pi} {\rm Area}(S)$
(one may also consult the argument of Lemma 3.1 in Marc Lackenby's paper \cite{L},
which can be easily modified to prove this inequality). A well-known computation
implies that ${\rm Area}(H)=l(\partial H)$. So we have

$$
\Sigma_i l_C(f(p_i)) \leq \Sigma_i l_H(p_i) = l(\partial H) =
$$ 
$$ {\rm Area}(H)\leq \frac{3}{\pi} 
{\rm Area}(S) = \frac{3}{\pi}\cdot 2\pi |\chi(S)| = 6|\chi(S)|$$
 
by the Gauss--Bonnet theorem. 
\end{proof}

\section{Word-hyperbolic Dehn filling}

Let $N$ be a finite volume hyperbolic 3--manifold with a unique embedded horocusp $C$. 
For a slope $\alpha$ in
$\partial C$, if $l_C(\alpha)>2\pi$, then Gromov and Thurston proved that 
$N(\alpha)$ has a metric of negative curvature.    
Theorem \ref{hyperbolike} implies that if $l_C(\alpha)> 6$, then $N(\alpha)$ is
hyperbolike. The intuition for why  such an
improvement is possible is that the $2\pi$--theorem only makes use of the negative
curvature of $N$ in the cusp $C$, whereas this result takes account of
negative curvature of $N$ outside of $C$
 as well. 

First, we need to state a theorem of Lackenby.
Let $N$, $C$, and $\alpha$ be as above. 
Let $k$ be the core of the Dehn filling in $N(\alpha)$. We will fix a Riemannian metric
on $N(\alpha)$ which agrees with the hyperbolic metric on $N\backslash C$. 
For a homotopically trivial mapping $c\co S^1\to N(\alpha)\backslash \mathcal N( k)$,
we define the {\it wrapping number} $$\text{wr}(c,k) = \min\{|d^{-1}\mathcal{N}(k)|,\ d\co D^2\to M(\alpha),
d\ \text{is transverse to $\mathcal{N}(k)$,}$$ $$\text{and}\ d_{|\partial D}=c\}.$$ 
It measures the minimal number of 
 intersections with $\mathcal{N}(k)$ of 
maps of disks spanning $c$. The following theorem is due to Lackenby \cite[Theorem 2.1]{L}:
\begin{theorem}[Ubiquity theorem]\label{ubiquity} 
In the situation above, 
there is a constant $w$ such that
for any least area disk $d\co D^2\to N(\alpha)$, we have $\text{area}(d)\leq 
w(\text{wr}(\partial D,k)
+\text{length}(\partial D))$. 
\end{theorem}

This theorem strengthens the ubiquity theorem of Gabai \cite{Ga2}, in that
it doesn't count the multiplicities of intersections of $d\co D^2\to N(\alpha)$
with $k$. The point of this theorem is that to obtain a linear isoperimetric
inequality for $N(\alpha)$, we need only show that there is a constant $v$ 
so that  for maps $d\co D^2\to N(\alpha)$, 
$\text{wr}(\partial D,k)\leq v \text{length}(\partial D))$. 

\begin{theorem}[Hyperbolike fillings] \label{hyperbolike}
Let $N$ be a finite volume hyperbolic 3--manifold with single embedded horocusp $C$. If 
\a\ is a slope with $l_C(\alpha) > 6$, then $N(\alpha)$ is hyperbolike.
\end{theorem}

\begin{proof}
If $\partial C$ is not embedded, replace $C$ with a slightly smaller horocusp, retaining
the property that $l_C(\alpha)>6$.
Define $M=N\backslash C$, and let $k$ be the core of the Dehn filling $N(\alpha)$,
and $\mathcal N(k)$ the open solid torus which is attached to $M$. 
Suppose $\pi_2 N(\alpha) \neq 0$ or $|\pi_1 N(\alpha)|<\infty$. $M$ has incompressible 
boundary, so by \ref{essential surfaces}, $N(\alpha)$ contains a mappping of a
sphere or disk $f\co S\to N(\alpha)$ 
such that $f_{|\hat{S}=f^{-1}(M)}$ is essential in $M$. 
Let $n=|f^{-1}(\mathcal N (k))|$. Then there are at least $n-1$ boundary components of $\hat S$ which 
map to  multiples of \a\ in \del$M$.
If $n = 0$, $\hat S$ would be an inessential sphere in $M$, since $\pi_2 M = 0$, and
therefore is trivial in $\pi_2 N(\alpha) $, a contradiction. If $n = 1$ or 2, so $\hat S$ is a 
disk or annulus, then $f_{|\hat S}$ can be homotoped 
into \del$M$, since \del$M$ is incompressible and $M$ is acylindrical. So $n\geq 3$. Applying Lemma 
\ref{bound}, we see
$$6(n-2)=6|\chi(S)|\geq (n-1)\cdot l_C(\alpha)>6(n-1)$$ a contradiction. So $N(\alpha)$ is
irreducible with the core having infinite order in $\pi_1N(\alpha)$. 

 Choose a metric on $N(\alpha)$
which agrees with the hyperbolic metric on $M$, 
and is any metric on $H=\overline{\mathcal{N}(k)}$.
We want to show that $N(\alpha)$ has linear isoperimetric inequality with this metric. 

Choose a map $c\co S^1 \to N(\alpha)$ which is homotopically trivial. 
First, we will find a homotopy of $c$ to a map $c'$ in $M$, such that the length $c'$
and the area of the homotopy
are linearly bounded by the length of $c$. The second step is to show that the wrapping number of 
$c'$ is linearly bounded by its length. We then apply the ubiquity theorem to conclude
that $N(\alpha)$ has linear isoperimetric inequality. 

Then $c^{-1}(\text{int} H)$ consists of a 
collection of intervals. 
Let us consider one of these intervals $\delta$. Lift $c_{|\delta}$ to a map
 $c_0\co \delta\to \tilde H$, where $\tilde H$ is the universal cover of $H$. Change
the metric on $H$ to be isometric to a euclidean cylinder quotient a 
translation. Then this
Riemannian metric is quasi-isometric to the original metric on $H$. 
We can homotope $c_0$ 
to a map $c_1\co \delta \to \partial\tilde H$ keeping endpoints fixed, such that 
$\text{ length}(c_1)\leq \frac{\pi}{2} \text{ length}(c_0)$
(see Figure  \ref{cylinder}), where $c_1(\delta)$ is a shortest arc in 
$\partial\tilde H$ connecting the endpoints of
$c_0(\delta)$.
\begin{figure}[ht!] 
\cl{\small
\SetLabels 
\E(.5*.6){$c_0(\delta)$}\\
\E(0.7*.73){$c_1(\delta)$}\\
\E(.4*.27){$\tilde H$}\\
\endSetLabels 
\AffixLabels{\includegraphics[width=3in]{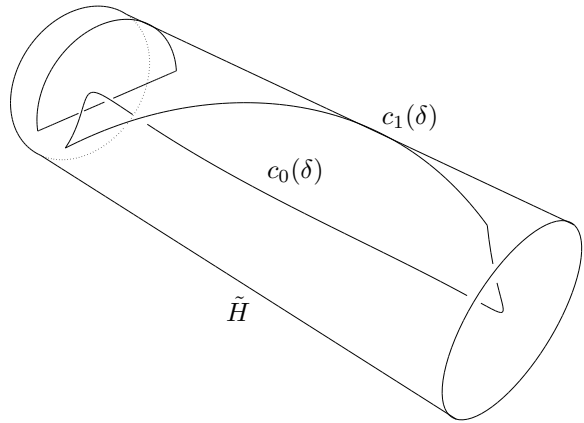}}}
	\caption{\label{cylinder} Comparing lengths}
\end{figure}
 The extremal case occurs when $c_0(\delta)$ is a diameter of the cylinder.
$c_0(\delta)\cup c_1(\delta)$ bounds a map of a disk whose area is linearly bounded by 
$\text{ length}(c_0)+\text{ length}(c_1)\leq C \text{ length}(c)$. For example,
the disk which connects each point of $c_0(\delta)\cup c_1(\delta)$
by the shortest segment to the axis of the cylinder works, as can
be seen by an elementary computation. 
Since the metric on $H$ is 
quasi-isometric to the euclidean metric,
we can find a homotopy of $c$ to a map of a curve  $c'\co S^1\to M$, 
whose length is linearly bounded by $c$, and 
such that the area of the homotopy is linearly bounded by $c$. Replace
$c$ with this map $c'$.

We want to 
estimate $\text{wr}(c,k)$, for $im(c)\subset M$. By lemma \ref{essential surfaces} we may assume 
that $c$ bounds  a map of a punctured disk $d\co S\to N$ such that $d_{|d^{-1}(M)}$ is incompressible and 
\del--incompressible in $M$, with $n$ boundary components of $d^{-1}(M)$ mapping to  
multiples of \a\ in $\partial M$, and $d_{|d^{-1}(C)}$ consists of maps of annuli which
can be assumed to be products with respect to the horotorus foliation of $C$. 
So $\text{wr}(c,k)\leq n$. If $c$ is homotopic to $\partial M$, then $c$ would be homotopic
in $M$ to a multiple of $\alpha$, since otherwise it would be homotopic to a multiple of $k$, and it would
not be homotopically trivial in $N(\alpha)$. The area of the annulus realizing the homotopy
into $\partial M$ can be chosen to be linearly bounded by $\text{ length}(c)$, for example
by coning off $c$ to the cusp in $N$. Therefore $c$ 
bounds a map of a disk in $N(\alpha)$  whose area is linearly bounded by $\text{ length}(c)$.
If $c$ is not homotopic into $\partial M$, then we may  homotope $c$ to be geodesic in $N$, and $d$ to 
be pleated in $N$, by lemma \ref{pleated surfaces}. 
Consider $d^{-1}(C)\subset S$. Then as in lemma
\ref{bound}, 
we can find disjoint cusp neighborhoods $H_i$ in $S$, some of which might intersect \del$S$.  
Let us estimate how many horocusps can
meet \del$S$.
We will assume the first $j$ cusps meet \del$S$. Shrink each cusp meeting \del$S$
 until it is tangent to \del$S$. Lifting to 
$\tilde S \subset \mathbb H ^2$, so that a geodesic component \g\ of $\widetilde{\partial S}$ 
runs from
$0$ to $\infty$, we see a sequence of $j+1$ horodisks tangent to \g, such that
the first and $j+1$st horodisks are identified by the covering translation of \g. So
the length of \del$S$ is the distance between the tangent points of these two horodisks. Consider two 
sequential horodisks. Then we may move the horodisk of larger euclidean radius by a hyperbolic
isometry, keeping it tangent 
to \g\ until  it is tangent to the smaller one. 
A geometric calculation shows that 

\begin{figure}[ht!]\cl{\small
\psfrag{A}{$R+r$}
\psfrag{B}{$R-r$}
\psfrag{C}{$R$}
\psfrag{D}{$R-r$}
\psfrag{E}{$r$}
\psfrag{F}{$\ln\frac{R}{r}$}
\psfrag{G}{}
\psfrag{H}{$2(R-r)^2=(R+r)^2$}
\psfrag{I}{$\Rightarrow\frac{R}{r}=(1+\sqrt{2})^2$}
\includegraphics[width=3.5in]{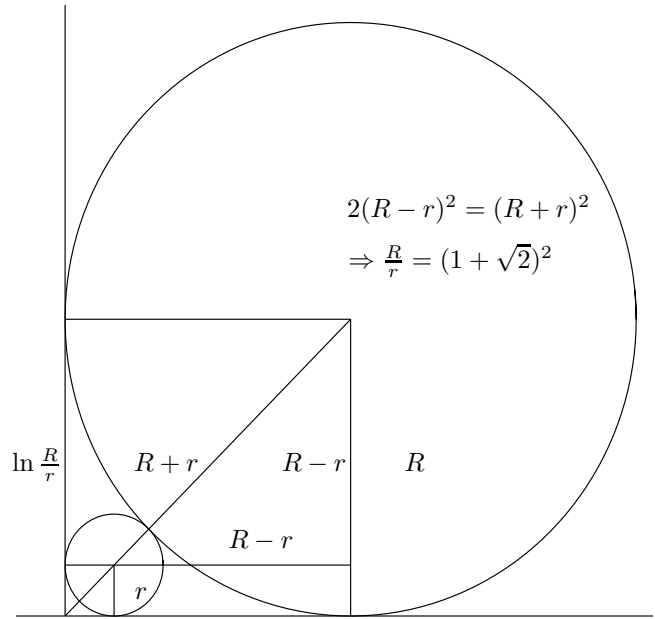}}
\caption{\label{horocalc} Bounding translation length}
\end{figure}
the hyperbolic length between the tangency points of the horodisks is $2\ln(1+\sqrt{2})$(Figure  \ref{horocalc}).
So $l(\partial S) \geq 2j\ln(1+\sqrt{2}) $. 

Take $S$ and double  it along its geodesic boundary $\partial S$  to a hyperbolic surface $DS$. 
As in lemma \ref{bound}, $l(\partial{H}_i)\geq l_C(\alpha)$, for $i> j$. So we take
the collection of horocusps in $DS$ consisting of $H_i$ and its
reflection, for $i> j$. Choose a number $\epsilon$ such that $l(\alpha)> 6+\epsilon$.
 Then we have
$$6(2n-2) = 6|\chi(DS)|\geq 2\sum_{i=j+1}^n l(\partial H_i) \geq 
2\sum_{i=j+1}^n l_C(\alpha)
\geq 2(n-j)(6+\epsilon).$$
Thus,
$$2\epsilon n\leq 2j(6+\epsilon)-12 \leq \frac{(6+\epsilon)l(c)}{\ln (1+\sqrt{2})}.$$
So $\text{ wr}(c,k)\leq n \leq \frac{(6+\epsilon)l(c)}{2\epsilon 
\ln (1+\sqrt{2})}$.
By the  ubiquity theorem \ref{ubiquity}, $N(\alpha)$ has linear 
isoperimetric inequality. 
\end{proof}

\section{Essential surfaces and Dehn filling}

The next theorem gives a condition for which a quasifuchsian surface
in a hyperbolic knot complement remains $\pi_1$--injective under Dehn filling. 

As usual, $N$ is a hyperbolic $3$--manifold with a horocusp $C$ and $S$ is a surface of finite type.  
Let $f\co S\to N $ be a $\pi_1$--injective mapping, taking
cusps of $S$ to cusps of $N$.
We will assume that 
the covering $\tilde N_f$ of $N$ corresponding to $f_*(\pi_1S)$ is geometrically finite. 
Let $Q(S)$ be the convex core of $\tilde N_f$.  
Suppose $f$ has no accidental parabolics, that is $Q(S)$ is homeomorphic
to $S\times [0,1]$ (or it is homeomorphic to $S$ if $\pi_1(S)$ is fuchsian), 
and all cusps of $S$ map to the same boundary slope in $C$ (could be none). 
Let $\tilde C$ be the preimage of $C$ in  $\tilde N_f$. 
Suppose $Q(S)\cap\ \tilde C \cong \mathcal{N}(cusps(Q(S)))$, 
that is the only intersections with $\tilde C$ are the ones which must occur. Call
such a mapping $f$ {\it geometrically proper with respect to C}.

\begin{theorem}[Quasifuchsian filling] \label{qf}
Assume we have $N$ and $f\co S\to N$ as above, so that $f$ is geometrically proper with
respect to $C$. 
Let \a\ be the  slope on $\partial C$ corresponding to the image of the cusps of $S$
under the mapping $f$, or any slope in $C$, if $Q(S)$ is compact. Suppose
 $l_C(\alpha)\geq 6$. Form the compact surface $S' \supseteq S\backslash f^{-1}(C)$ such
that $K=S'\backslash( S\backslash f^{-1}(C))$ consists of disks, and a mapping $f'\co S'\to N(\alpha)$,
such that $f'_{|S'\backslash K}=f$ and $f'_{|K}\subset
N(\alpha)\backslash (N\backslash C)$. Then $f'$ is $\pi_1$--injective in $N(\alpha)$. 
\end{theorem}
\begin{proof}
Suppose $f'$ is not injective into  $\pi_1 N(\alpha)$. Let  $g\co S^1\to S'\backslash K$ be a
map which is homotopically non-trivial in $S'$
and which 
bounds a map of a disk $D$ into $N(\alpha)$, that is there is a map $d\co D\to N(\alpha)$
with $d_{|\partial D}=f\circ g$.  Choose $d^{-1}(\mathcal N(k))$ to have as few components 
as possible, where $k$ is the core of the Dehn filling.  Then $f\circ g$ is homotopic 
to a unique map with geodesic
 image \g\ in $N$,
 which will lie
inside of $Q(S)$ when we lift to $\tilde{N_f}$, since $g$ is homotopically non-trivial in $S'$. 
By lemmas \ref{essential disks} and \ref{pleated surfaces} \g\ bounds an incompressible, 
\del--incompressible pleated map of a punctured disk
$d\co  F \to N$, with $n$ punctures mapping to  multiples of \a\ in $C$, $d(\partial F)=\gamma$. 
Suppose 
$d^{-1}(C)\cap \partial F \neq \emptyset$. Then look at a component of $d^{-1}(C)$ which
intersects $\partial F$, and suppose it is noncompact. Then there is an embedded arc $\beta$ in  $d^{-1}(C)$
connecting a point in $\partial F$ to a cusp of $F$. There must also be such a geodesic arc $\beta'$ in
$Q(S)\cap \tilde C$ in the cover $\tilde N_f$ connecting the preimage of \g\ with the corresponding 
cusp in $\tilde C$, by the assumption that $f$ is geometrically proper with respect to $C$. 
Since the lift of $d(\beta)$ to $\tilde N_f$ and $\beta'$ lie entirely in
the same component of $\tilde C$, $d$ can be homotoped so that $d(\beta)=\beta'$. 
Take a neighborhood $R$ of $\beta$ in $F$ which contains the cusp at one end of 
$\beta$, and consider the subsurface $F'=F\backslash R$. Then $d_{|\partial F'}$
lifts to a map into $Q(S)$, and there is a map $g'\co S^1\to S'\backslash K$ 
so that $f\circ g'$ is homotopic to $d_{|\partial F'}$. Moreover, $g'$ is
homotopic to $g$, since $\pi_1 Q(S)=\pi_1(S)$. Thus, we have found a 
map of a loop $g'\to S'\backslash K$ which bounds a map of a disk in $N(\alpha)$
with fewer intersections with $\mathcal N(k)$
Thus, every component of $d^{-1}(C)$ which intersects $\partial F$ must be compact. If $n = 0$ or 1, 
then $S$ would be compressible, or have an accidental parabolic. Otherwise we may apply the argument
of theorem
\ref{bound} to get
$$6(n-1) = 6|\chi(F)|\geq n\cdot l_C(\alpha)\geq 6n$$
a contradiction. The point is that since $d^{-1}(C)$ intersects $\partial F$ only in
compact pieces, we may find embedded horocusp neighborhoods of the punctures in $F$
which miss $\partial F$, so that we may apply Boroczky's theorem to the double
of $F$, as we did in theorem \ref{hyperbolike}. 
\end{proof}

\begin{figure}[ht!] 
	\begin{center}
	\includegraphics[width=2.4in]{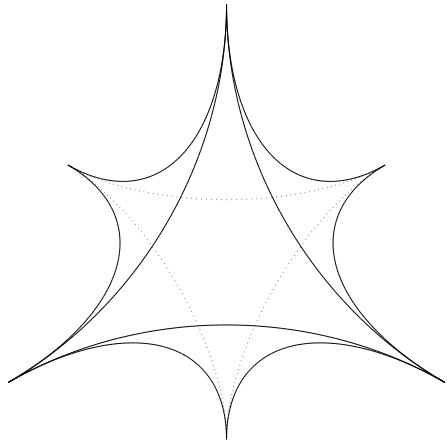}
	\caption{\label{octahedron} An ideal octahedron O in the conformal model}
	\end{center}
\end{figure}

Here is an example which shows that the bound in theorem \ref{hyperbolike} is 
sharp. We construct a manifold which has a totally geodesic punctured torus
with maximal possible cusp size. Take an ideal octahedron $O$ in $\H$ which has all angles
between faces $\pi/2$, as in Figure  \ref{octahedron}.

\begin{figure}[ht!] 
	\begin{center}\small
	\psfrag{A}{A}
	\psfrag{B}{B}
	\psfrag{C}{C}
	\psfrag{a}{a}
	\psfrag{b}{b}
	\psfrag{c}{c}
	\psfrag{d}{d}
	\psfrag{e}{e}
	\psfrag{f}{f}
	\includegraphics[width=3.6in]{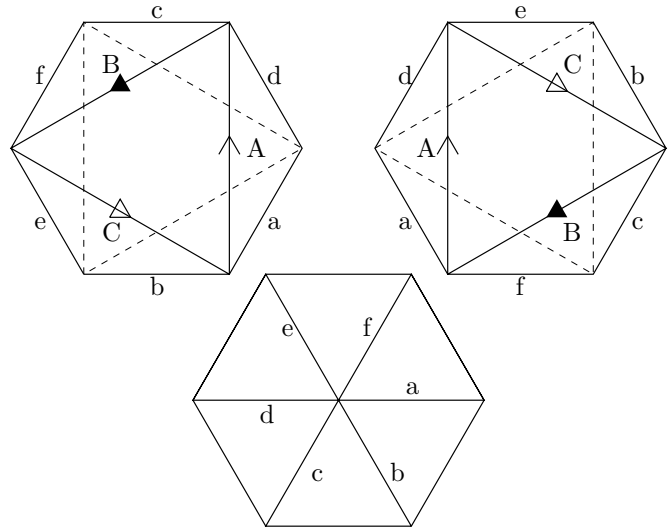}
	\caption{\label{manifold} How to glue up the manifold}
	\end{center}
\end{figure}

 Then we take two copies of $O$, and glue the top six side faces
 together in pairs
as indicated in Figure  \ref{manifold}.
 The edges of the front faces get glued up in such a way that we get a 
punctured torus made of two ideal triangles. The six side
faces get glued cyclically,
to form a punctured disk, as in the bottom diagram of figure \ref{manifold}. 
Double the manifold obtained so far along this punctured disk, then there are
two punctured tori, and the back faces of the two octahedra double to form two 3--punctured spheres. 
We then glue the punctured tori and punctured spheres together to get a manifold $N$ of finite 
volume, with 4 cusps (we can get two cusps by gluing the punctured spheres
with a twist). The cusp $C$ corresponding to the punctured torus has an embedded horoball
neighborhood with boundary slope length = 6. The punctured torus remains incompressible
after Dehn filling along this slope, by theorem \ref{qf} (this can also be
shown using the fact that the torus is homologically non-trivial, and the
filling is irreducible).  This shows that the bound given in \ref{bound} is sharp. 
By Dehn filling the other cusps of $N$, we can get manifolds with an embedded
punctured torus and a cusp corresponding to $C$, such that the boundary slope
is as close to $6$ as we like. This shows that the theorem \ref{hyperbolike} is
sharp as well.  

\begin{figure}[ht!] 
\cl{\small
\SetLabels 
\E(.37*.5){$\frac1n$}\\
\E(.6*.05){$\frac1n$}\\
\E(.6*.95){$-\frac1n$}\\
\E(.96*.5){$-\frac1n$}\\
\endSetLabels 
\AffixLabels{\includegraphics[width=3.5in]{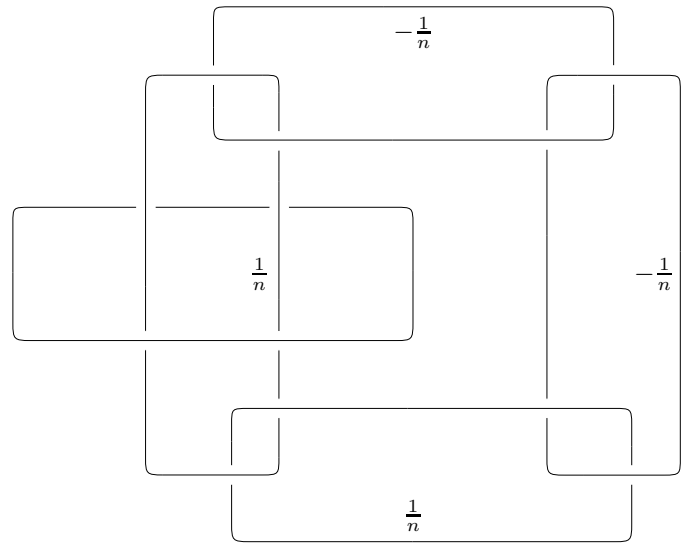}}}
\caption{\label{waistknot} Knots with meridian length $\to 4$ as 
$|n|\to \infty$}
\end{figure}

Here is an example of hyperbolic knots in $S^3$ with meridian slope length
in a maximal horocusp approaching 4. Take the 5 component link $L$ which is the 2--fold
branch cover over one component of the Borromean rings. It is well known that
the meridian slope for a maximal horocusp in the Borromean rings is 2, so the 
link $L$ has one component with meridian slope length 4. Then we may do
arbitrarily high Dehn fillings on the other components to obtain knots in 
$S^3$ with meridian slope length approaching 4 (see figure \ref{waistknot} for
the Dehn filling description). One may see that the Dehn fillings in diagram \ref{waistknot} 
on each pair of unlinked components cancel each other by opposite Dehn twists
on an annulus connecting up each pair in the complement of the other pair, 
so that the manifold obtained by the Dehn filling is still a knot 
in $S^3$.  It would be interesting to find knots with
longer meridian slope lengths. 

\section{Bounds on exceptional slopes}

For a pair of slopes \a, \b\ on a torus, call their intersection number $\Delta(\alpha,\beta)$.
If we choose a basis for the homology on the torus, such that $\alpha=(a,b),\ \beta=(c,d),$
then $\Delta(\alpha,\beta)=|ad-bc|$. If $\alpha$ is a slope, then $\gcd(a,b)=1$, since
$\alpha$ represents a primitive homology class. 

\begin{theorem}
Let $N$ be a hyperbolic $3$--manifold, and $C$ a distinguished embedded torus cusp. 
The intersection number  be\-tween exceptional boundary slopes on $C$ is $\leq 10$, and
there are at most 12 exceptional boundary slopes. 
\end{theorem}
\begin{proof}
Given two exceptional slopes \a, $\beta$, $l_C(\alpha)\leq 6$, and $l_C(\beta)\leq 6$ by
theorem \ref{hyperbolike}.
By a result of Cao and Meyerhoff, theorem 5.9 in \cite{CM}, $\text{area}(\partial C)\geq 3.35$. 
Let $\theta_{\alpha\beta}$ be the angle between the geodesics \a\ and $\beta$ on \del$C$. 
Computing area, we have $l_C(\alpha)\cdot l_C(\beta) \sin(\theta_{\alpha\beta}) = \Delta (\alpha,\beta)
\text{area}(\partial C)$. So $$\Delta(\alpha,\beta) = \frac{l_C(\alpha)\cdot l_C(\beta) 
\sin(\theta_{\alpha\beta})}{\text{area}(\partial C)} \leq \frac{6^2}{3.35}=10.75,$$ so 
$\Delta(\alpha,\beta)\leq 10$.

For the second part of the claim, we need the following lemma:
\begin{lemma}[Bound on number of slopes] 
If a collection of slopes on a torus have pairwise  intersection numbers
$\leq R$, then for  any prime number $p > R$, the number of such slopes is
bounded by $p+1$.
\end{lemma}
\begin{proof}
 Denote the projective plane over the finite 
field of order $p$ by
$\fp$. Then there is a map 
$\qp\to\fp$, where 
$\frac{a}{b}\mapsto (a\ \mod\ p,b\ \mod\ p)$. This map is well-defined, since 
if $\frac{a}{b}\mapsto (0,0)$, 
then $p|\gcd(a,b)=1$. Suppose a pair of slopes  
$\frac{a}{b}$ and $\frac{c}{d}$ in the given collection map to the 
same point in $\fp$, then $(a,b)\equiv k(c,d) (\mod\ p)$, so $|ad-bc|\equiv |kcd-kcd| \equiv 0\ (\mod\ p)$. 
If $|ad-bc|=0$, then  
$\frac{a}{b}=\frac{c}{d}$. Otherwise, 
$$p\leq |ad-bc|= \Delta(\frac{a}{b},\frac{c}{d})\leq R < p,$$
a contradiction. So for each point of $\fp$, there is at most one slope in the
collection mapped to it. Thus, there are at most $|\fp|=p+1$ slopes
in the collection. 
\end{proof}

In the case at hand, we have $R=10<11$, so we compute that the 
number of exceptional fillings is
 $\leq 12$. 
\end{proof}

\begin{figure}[ht!] 
	\begin{center}
	\includegraphics{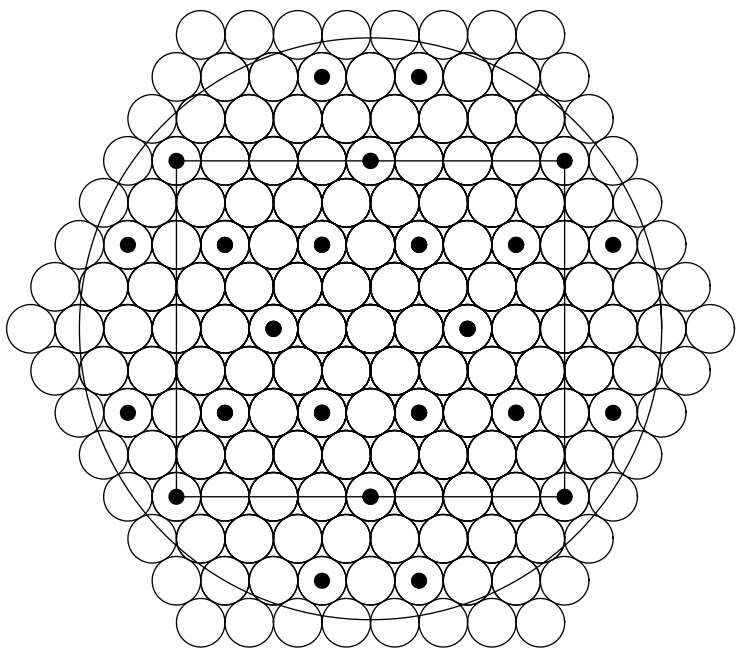}
	\caption{\label{8sister} Primitive lattice points in the figure eight knot sister}
	\end{center}
\end{figure}

It is conjectured that the maximal intersection number between 
exceptional slopes is 8, realized by the figure 8 knot complement \cite{G}.
Moreover, we expect that the figure 8 knot has the fewest number 
of exceptional slopes, 10. When applied to the figure 8 knot, theorem
\ref{hyperbolike} gives exactly the set of exceptional slopes for
the maximal cusp.
On the other hand, the figure eight knot 
sister has a regular torus cusp,
with 12 slopes of length $\leq 6$, but there are only 8 exceptional 
fillings \cite{G} (see Figure  \ref{8sister}). In the figure, the view is 
from $\infty$ in the cusp of the figure eight sister, and the circles
correspond to other horoball copies of the cusp from our viewpoint.  Signed pairs of lattice points
correspond to slopes, where a segment from the center of the picture to the
lattice point maps down to a boundary slope in the manifold. The exceptional slopes are shown
in the box. 

\end{document}